%% file: RoysetByun_revised.tex
\documentclass[titlepage,final,10pt]{article}
\input mac.tex

\begin{document}


\begin{center}
\begin{large}
{\bf Gradients and Subgradients of Buffered Failure Probability}
\smallskip
\end{large}
\vglue 0.7truecm
\begin{tabular}{cc}
  \begin{large} {\sl Johannes O. Royset} \end{large} & \begin{large} {\sl Ji-Eun Byun } \end{large} \\
  Operations Research Department ~~&~~    Engineering Risk Analysis Group\\
  Naval Postgraduate School &    Technical University of Munich\\
  joroyset@nps.edu & j.byun@tum.de
\end{tabular}

\vskip 0.2truecm

\end{center}

\vskip 1.3truecm

\noindent {\bf Abstract}. \quad Gradients and subgradients are central to optimization and sensitivity analysis of buffered failure probabilities. We furnish a characterization of subgradients based on subdifferential calculus in the case of finite probability distributions and, under additional assumptions, also a gradient expression for general distributions. Several examples illustrate the application of the results, especially in the context of optimality conditions.\\

\vskip 0.5truecm

\halign{&\vtop{\parindent=0pt
   \hangindent2.5em\strut#\strut}\cr
                         \cr\cr

{\bf Date}:\quad \ \today \cr}

\baselineskip=15pt

\section{Introduction}

In decision making under uncertainty, it is common to assess a random quantity of interest by computing the probability of outcomes above a threshold. With little loss of generality, one can shift the quantity of interest and adopt zero as the threshold. This leads to the {\em probability of failure} with numerous uses in reliability engineering and stochastic optimization; see \cite{KucukyavuzRuiwei.21,RockafellarRoyset.15b} and references therein. Despite advances in the analysis and computation of gradients and subgradients of failure probabilities with respect to parameters \cite{HantouteHenrionPerezaros.19,PenaordieresLuedtkeWachter.20,RoysetPolak.07,Uryasev.95,vanackooijJavalPerezaros.21}, it remains theoretically and computationally challenging to use failure probabilities in sensitivity analysis and optimization. Convexity of functions describing quantities of interest does not necessarily translate into convexity of failure probabilities \cite{RockafellarRoyset.10}. In a data-driven setting with finite probability distributions, a failure probability is expressed by a weighted average of indicator functions and thus continuity is lost as well. One tends to end up with gradients that are either zero or not defined, making sensitivity analysis and optimization difficult. These challenges motivate the alternative {\em buffered failure probability} \cite{RockafellarRoyset.10}. Analogous to the shift from quantiles to superquantiles pioneered in \cite{RockafellarUryasev.00}, the pass from conventional failure probabilities to buffered failure probabilities results in several desirable properties \cite{MafusalovUryasev.18}. In this paper, we provide gradient and subgradient expressions for buffered failure probabilities and discuss their use in sensitivity analysis and optimization.

Buffered failure probabilities and the closely related buffered probabilities of exceedance \cite{DavisUryasev.16} appear in several optimization models, often leveraging their connection with superquantiles (a.k.a. conditional/average value-at-risk), with applications in reliability engineering \cite{BasovaRockafellarRoyset.11,ByunRoyset.21,ChaudhuriKramerNortonRoysetWillcox.21,RockafellarRoyset.10,ZrazhevskyGolodnikovUryasevZrazhevsky.20}, data analytics \cite{NortonUryasev.19}, finance \cite{ShangKuzmenkoUryasev.18} and to address natural disasters \cite{DavisUryasev.16}. However, gradient expressions for buffered failure probabilities as functions of parameters are not fully developed. Since a buffered failure probability can be viewed as the ``inverse'' of superquantiles, one can leverage gradient expressions for the latter and the implicit function theorem to compute gradients of buffered failure probabilities \cite{ZhangUryasevGuan.19}. While this suffices in some cases, a superquantile is not generally smooth and sufficient conditions for such smoothness are nontrivial. Under a convexity assumption about the quantity of interest, subgradients of superquantiles emerge as quasigradients of buffered failure probabilities \cite{ZhangUryasevGuan.19}. However, engineering systems often involve nonconvex quantities of interest; see for example \cite{BasovaRockafellarRoyset.11}. When a finite distribution produces distinct outcomes for the quantity of interest and certain expectations are bounded away from zero, gradient expressions for buffered failure probabilities become available \cite{ByunRoyset.21}. We extend on these developments by considering both finite and continuous distributions, and circumventing any assumption about convexity. Our proofs of gradients and subgradients avoid relying on the implicit function theorem and distinct outcomes and instead leverage subdifferential calculus as laid out in \cite{VaAn,primer} and properties of expectation functions in \cite{ShapiroDentchevaRuszczynski.09}.

For a quantity of interest expressed by $g:\reals^m\times\reals^n\to \reals$ and an $m$-dimensional random vector $\bfxi$, the (conventional) failure probability under decision $x\in\reals^n$ is given by
\[
p(x) = \prob\big\{ g(\bfxi,x) > 0 \big\},
\]
which is well defined as long as $g(\,\cdot\,, x)$ is measurable. (We use boldface letters to indicate a random quantity and return to regular font for its outcomes.) Under the assumption that $g(\bfxi,x)$ is integrable (which we retain throughout), the {\em buffered failure probability} is defined as
\[
\bar p(x) = \bprob\big\{ g(\bfxi,x) > 0 \big\} = \begin{cases}
0 & \mbox{ if } p(x) = 0\\
1-\bar \alpha & \mbox{ if } p(x) > 0, ~~ \Ex\big[g(\bfxi,x)\big] < 0\\
1 & \mbox{ otherwise,}
\end{cases}
\]
where $\bar \alpha$ is the probability level $\alpha\in [0,1]$ that makes the $\alpha$-superquantile of $g(\bfxi,x)$ equal to zero. Although not explicitly indicated, $\bar\alpha$ typically depends on $x$. Recall that the $\alpha$-superquantile of $g(\bfxi,x)$ is given by
\[
\bar q_\alpha(x) = q_\alpha(x) + \frac{1}{1-\alpha}\Ex\Big[\max\big\{0, g(\bfxi,x) - q_\alpha(x) \big\}\Big]
\]
for $\alpha\in [0,1)$, where $q_\alpha(x)$ is the $\alpha$-quantile of $g(\bfxi,x)$, and $\bar q_1(x) = \lim_{\alpha\upto 1} \bar q_\alpha(x)$; see \cite[Section 3.C]{primer}.

As a key departure from \cite{ZhangUryasevGuan.19}, we leverage an alternative formula due to Norton and Uryasev \cite{NortonUryasev.19}:
\begin{equation}\label{eqn:minformula}
\bar p(x) = \nmin_{\gamma\geq 0} \,\Ex\Big[\max\big\{0, ~\gamma g(\bfxi,x) + 1 \big\}\Big],
\end{equation}
which applies under the assumption that $p(x)>0$; see \cite[Proposition 2.2]{MafusalovUryasev.18}. The set of minimizers,
\begin{equation}\label{eqn:gamma}
\Gamma(x) = \nargmin_{\gamma\geq 0} \,\Ex\Big[\max\big\{0, ~\gamma g(\bfxi,x) + 1 \big\}\Big],
\end{equation}
contains $-1/q_{\bar \alpha}(x)$ when $\Ex[g(\bfxi,x)] < 0$ and $p(x) > 0$, where $\bar\alpha$ is such that $\bar q_{\bar\alpha}(x) = 0$; see \cite[Proposition 1]{NortonUryasev.19}. In fact, $-1/q_{\bar \alpha}(x)$ is the only minimizer unless the distribution function $P$ of $g(\bfxi,x)$ has a flat stretch to the right of $q_{\bar \alpha}(x)$ and $P(q_{\bar \alpha}(x)) = \bar \alpha$ as explained in the same reference.

The paper is laid out as follows. Section 2 addresses the case of finite probability distributions. Section 3 presents sufficient conditions for the buffered failure probability to be smooth and gives an expression for the gradients. In both cases, we assume that $g(\xi,\cdot\,)$ is smooth, at least locally near a point $x$ of interest but not necessarily for all $\xi$ as specified in detail below. Convexity is not required anywhere in the paper. These assumptions are well aligned with applications in engineering settings \cite{RockafellarRoyset.10,RockafellarRoyset.15b}. The paper ends in Section 4 with application of the results, especially to derive optimality conditions.

\section{Subgradients under Finite Distributions}

When the distribution of $\bfxi$ is finite with outcomes $\xi^1, \dots, \xi^\nu$, each occurring with a positive probability, the buffered failure probability $\bar p$ is differentiable at $\hat x$ under the assumptions that $\{g(\xi^i,\hat x), i=1, \dots, \nu\}$ are all distinct, certain expectations are bounded away from zero and $g(\xi^i,\cdot\,)$ is smooth (i.e., continuously differentiable) for all $i$ in a neighborhood of $\hat x$; see Example 1 and \cite{ByunRoyset.21}. Here, we derive a necessary condition for subgradients of $\bar p$ at $\hat x$ without requiring distinct values of $\{g(\xi^i,\hat x), i=1, \dots, \nu\}$ and the expectation assumptions. The set of subgradients $\partial f(x)$ of a function $f:\reals^n\to [-\infty, \infty]$ is well defined at any point $x$ with $f(x)$ being finite. We note that subgradients are understood to be of the general kind \cite[Chapter 8]{VaAn}, also referred to as limiting or Mordukhovich subgradients; see \cite[Section 4.I]{primer} for an introduction.

\begin{theorem}\label{thm:discrete}
    Suppose that a point $\hat x\in \reals^n$, an $m$-dimensional random vector $\bfxi$ and a function $g:\reals^m\times\reals^n\to \reals$ satisfy the following properties:
  \begin{enumerate}[{\rm (a)}]

  \item $\bfxi$ has a finite distribution with outcomes $\xi^1, \dots, \xi^\nu$ and corresponding positive probabilities $p_1, \dots, p_\nu$.

  \item For all $i=1, \dots, \nu$, the function $g(\xi^i,\cdot\,)$ is smooth in a neighborhood of $\hat x$.

  \item $\Ex[g(\bfxi,\hat x)] < 0$ and $g(\xi^i,\hat x) > 0$ for some $i$.

  \end{enumerate}
If $y \in \partial \bar p(\hat x)$, then there are $\hat\gamma\in \Gamma(\hat x)$ and multipliers $\mu_i \in [0,1], i=1, \dots, \nu$, such that
\[
y = \hat\gamma \sum_{i=1}^\nu p_i \mu_i \nabla_x g(\xi^i,\hat x).
\]
The multipliers satisfy the properties
\begin{equation}\label{eqn:multiplercond}
\sum_{i=1}^\nu p_i \mu_i g(\xi^i,\hat x) = 0~~~~~~~\mbox{ and,~~ for all $i$: } ~~~~~ \mu_i \in  \begin{cases}
  \{0\} & \mbox{ if } \hat\gamma g(\xi^i, \hat x) + 1 < 0\\
  [0,1] & \mbox{ if } \hat\gamma g(\xi^i, \hat x) + 1 = 0\\
  \{1\} & \mbox{ if } \hat\gamma g(\xi^i, \hat x) + 1 > 0.
  \end{cases}
  \end{equation}
\end{theorem}
\state Proof. Let $i_*$ be such that $g(\xi^{i_*},\hat x)>0$, which exists in light of assumption (c). By assumption (b), there is a nonempty compact set $X\subset\reals^n$ such that $\hat x$ is in its interior, $g(\xi^i, \cdot\,)$ is smooth at every $x\in X$ for all $i$ and $g(\xi^{i_*},x) \geq \half g(\xi^{i_*},\hat x)$ for all $x\in X$. We consider the function $f:\reals^n\times \reals\to (-\infty, \infty]$ given by
\[
f(x,\gamma) = \begin{cases}
\sum_{i=1}^\nu p_i f_i(x,\gamma) & \mbox{ if } x\in X, \gamma \geq 0\\
\infty & \mbox{ otherwise,}
\end{cases}
\]
where $f_i(x,\gamma) = h(\phi_i(x,\gamma))$, $\phi_i(x,\gamma) = \gamma g(\xi^i,x) + 1$ and $h(\eta) = \max\{0,\eta\}$. Certainly,
\[
\bar p(x) = \ninf_{\gamma \in \reals} f(x,\gamma) ~~~\forall x\in X.
\]
We now bring in subdifferential calculus rules for such inf-projections. Let $y \in \partial \bar p(\hat x)$.  By \cite[Theorem 5.13]{primer}, there exists a minimizer $\hat\gamma$ of $f(\hat x,\cdot\,)$ such that $(y,0) \in \partial f(\hat x,\hat\gamma)$. It is clear from the construction of $f$ that $\hat\gamma\in \Gamma(\hat x)$. The cited theorem requires that $f$ is proper and lsc, which is obvious in view of (b) and the construction of $X$. The theorem also requires that for all $\beta\in \reals$ and $\bar x\in \reals^n$ there are $\epsilon >0$ and a bounded set $C \subset \reals$ such that the level-sets $\{\gamma~|~f(x,\gamma) \leq \beta\} \subset C$ whenever $\|x-\bar x\|_2 \leq \epsilon$. By construction of $X$,
\[
f(x,\gamma) \geq p_{i_*} \big(\gamma g(\xi^{i_*},x) + 1 \big) \geq p_{i_*} \big(\half \gamma g(\xi^{i_*},\hat x) + 1 \big) ~~~\forall x\in X, ~\gamma\geq 0.
\]
Thus, we can take $C = [0, ~2\beta/(p_{i_*}g(\xi^{i_*},\hat x))]$ and the theorem applies. It remains to compute the subgradients of $f$ at $(\hat x, \hat \gamma)$. We note that $\hat \gamma >0$ because if $\hat\gamma$ were zero, then $\bar p(\hat x) = 1$ by \eqref{eqn:minformula}. This in turn entails that $\Ex[g(\bfxi,\hat x)] \geq 0$ by the definition of the buffered failure probability. However, this is ruled out by assumption (c).

First, we consider the subgradients of $f_i$. Without loss of generality, we can assume that $\phi_i$ is smooth because it can be extended beyond $X\times [0,\infty)$ in a smooth manner. Thus, the chain rule \cite[Theorem 4.64]{primer} applies and we obtain that
\[
\partial f_i(\hat x,\hat\gamma) = \nabla \phi_i(\hat x,\hat \gamma) \partial h\big(\phi_i(\hat x,\hat \gamma)\big) = \begin{cases}
  \{(0,0)\} & \mbox{ if } \hat \gamma g(\xi^i,\hat x) + 1 < 0\\
\big(\hat\gamma \nabla_x g(\xi^i,\hat x), ~g(\xi^i,\hat x) \big)  [0,1] & \mbox{ if } \hat\gamma g(\xi^i,\hat x) + 1 = 0\\
  \big\{\big(\hat \gamma \nabla_x g(\xi^i,\hat x), ~g(\xi^i,\hat x) \big)\big\} & \mbox{ if } \hat \gamma g(\xi^i,\hat x) + 1 > 0.
  \end{cases}
\]
Second, by the sum rule \cite[Proposition 4.67]{primer}, we achieve the expression
\[
\partial f(\hat x,\hat \gamma) = \sum_{i=1}^\nu p_i \partial f_i(\hat x, \hat\gamma),
\]
which holds by an equality, and not merely an inclusion, because $f_i$ is epi-regular at $(\hat x, \hat\gamma)$. This means that
\[
(y,0) = \sum_{i=1}^\nu p_i \mu_i \big(\hat \gamma \nabla_x g(\xi^i,\hat x), ~g(\xi^i,\hat x) \big),
\]
where $\mu_i$ is specified in the theorem. The conclusion follows.\eop

Beyond the restriction to finite distributions, which is common in a data driven setting, the theorem is widely applicable. Assumption (c) rules out pathological cases and is hardly a limitation.\\

\state Example 1. Consider a Bernoulli random variable $\bfxi$ with outcomes $\{0, 1\}$ and equal probabilities. Let $g(\xi,x) = x\xi  -1$ and consider $\hat x \in (1,2)$. This means that $\Ex[g(\bfxi,\hat x)] = \hat x/2 - 1 < 0$ and $g(1,\hat x)= \hat x - 1 > 0$. Direct computation establishes that the $\alpha$-quantiles and  $\alpha$-superquantiles are given by
\[
q_\alpha(\hat x) = \begin{cases}
  -1 & \mbox{ if } \alpha \in (0,1/2]\\
  \hat x-1 & \mbox{ if } \alpha \in (1/2,1)
\end{cases}
~~~~~~~~
\bar q_\alpha(\hat x) = \begin{cases}
  -1 + \hat x/(2(1-\alpha)) & \mbox{ if } \alpha \in [0,1/2]\\
  \hat x-1 & \mbox{ if } \alpha \in (1/2,1].
\end{cases}
\]
We solve the equation $\bar q_\alpha(\hat x) = 0$ and find $\bar \alpha = 1-\hat x/2$ so that $\bar p(\hat x) = \hat x/2$  and $\nabla \bar p(\hat x) = 1/2$.

Next, let's examine Theorem \ref{thm:discrete} in this case. We find that $y\in \partial \bar p(\hat x)$ implies that there are a minimizer $\hat\gamma\in \Gamma(\hat x)$ and multipliers $\mu_1$ and $\mu_2$
such that
\begin{align*}
  y & = \frac{1}{2} \hat\gamma \cdot 0 \cdot \mu_1 + \frac{1}{2} \hat\gamma \cdot 1 \cdot \mu_2  =\frac{1}{2}\hat \gamma\mu_2\\
  0 & = \frac{1}{2} (-1) \mu_1 + \frac{1}{2} (\hat x-1) \mu_2.
\end{align*}
From \eqref{eqn:gamma} and the surrounding discussion, $\hat\gamma = -1/q_{\bar\alpha}(\hat x) = 1$. This minimizer is unique because $\bar \alpha \in (0,1/2)$, which, in turn, implies that the distribution function of $g(\bfxi,\hat x)$ is strictly higher than $\bar \alpha$ at $q_{\bar\alpha}(\hat x)$. The requirement from the theorem about $\mu_1$ and $\mu_2$ translates into $\mu_1\in [0,1]$ and $\mu_2 = 1$ because
\begin{align*}
  \hat\gamma (\hat x\cdot 0 - 1) + 1 & = 0\\
  \hat\gamma (\hat x\cdot 1 - 1) + 1 & = \hat x > 0.
\end{align*}
Thus, $y = 1/2$ so there is only one subgradient and this is indeed the gradient $\nabla \bar p(x)$.

For comparison, let us consider the gradient formula in \cite{ByunRoyset.21}:
\begin{equation}\label{eqn:firstpaper}
\nabla \bar p(x) =\sum_{i=j+1}^\nu p_i \Bigg(\frac{g(\xi^i,x)}{(g(\xi^{j},x))^2}\nabla_x g(\xi^{j},x) -\frac{1}{g(\xi^{j},x)}\nabla_x g(\xi^i, x)\Bigg),
\end{equation}
where $g(\xi^1,x) < g(\xi^2,x) < \cdots < g(\xi^\nu,x)$ are sorted and $j \in \{1, \dots, \nu\}$ satisfies
\[
\sum_{i=j}^\nu p_i g(\xi^i,x) < 0 < \sum_{i=j+1}^\nu p_i g(\xi^i,x).
\]

For the present example, $g(0,\hat x) = -1$ and $g(1,\hat x) = \hat x-1$. Thus, we can set $j = 1$ because
\[
\tfrac{1}{2} g(0,\hat x) + \tfrac{1}{2} g(1,\hat x) = \tfrac{1}{2} \hat x - 1 < 0 < \tfrac{1}{2} g(1,\hat x) = \tfrac{1}{2} \hat x - \tfrac{1}{2}.
\]
Using the formula \eqref{eqn:firstpaper}, we obtain
\[
\nabla \bar p(\hat x) = \frac{1}{2}\Bigg(\frac{g(1,\hat x)}{(g(0,\hat x))^2}\nabla_x g(0,\hat x) - \frac{1}{g(0,\hat x)}\nabla_x g(1, \hat x)  \Bigg) = \frac{1}{2}
\]
as expected. Example 3 furnishes a case when \eqref{eqn:firstpaper} does not apply.\eop

\section{Gradients of Buffered Failure Probability}

Since a buffered failure probability is expressed in terms of ``min'' and ``max'' (see \eqref{eqn:minformula}), we cannot expect it to define a smooth  function in general. However, the following requirements on the function $g$ and the distribution of $\bfxi$ ensure that $\bar p$ is differentiable at a particular point $\hat x$.

\begin{theorem}\label{thm:contcase}
  For a point $\hat x\in \reals^n$, an $m$-dimensional random vector $\bfxi$ and a function $g:\reals^m\times\reals^n\to \reals$, suppose that there is a neighborhood $X$ of $\hat x$ such that the following hold:
  \begin{enumerate}[{\rm (a)}]

  \item For each $x\in X$, there is a set $\Xi\subset\reals^m$, with $\prob\{\bfxi\in \Xi\} = 1$, such that $g(\xi, \cdot\,)$ is smooth at $x$ for every $\xi\in \Xi$.

  \item There are $\epsilon>0$ and $\delta\in (0,1]$ such that $\prob\{g(\bfxi,x)\geq \epsilon\} \geq \delta$ for all $x\in X$.

  \item There is a measurable function $h:\reals^m\to [0,\infty)$ such that $\Ex[h(\bfxi)] < \infty$ and $|g(\xi,x)| \leq h(\xi)$ for all $\xi\in\reals^m$ and $x\in X$.

  \item There is a measurable function $\kappa:\reals^m\to [0,\infty)$ such that $\Ex[\kappa(\bfxi)] < \infty$ and
   \[
   \big| g(\xi,x) - g(\xi,x') \big| \leq \kappa(\xi) \|x-x'\|_2   ~~~\forall x,x'\in X, ~~\xi\in \reals^m.
   \]

  \item If $\gamma>0$ and $x\in X$, then $\prob\{-\gamma g(\bfxi,x) = 1\} = 0$.

  \item The set of minimizers $\Gamma(\hat x)$ is the single point $\hat \gamma$.
  \end{enumerate}
Then, $\bar p$ is differentiable at $\hat x$ and
  \[
  \nabla \bar p(\hat x) = \Ex\big[ F(\bfxi,\hat x,\hat\gamma) \big], ~~\mbox{ where } F(\xi,\hat x,\hat\gamma) = \begin{cases}
  0 & \mbox{ if } \hat\gamma g(\xi,\hat x) + 1 \leq 0\\
  \hat\gamma \nabla_x g(\xi,\hat x) & \mbox{ otherwise}.
  \end{cases}
  \]
\end{theorem}
\state Proof. Assumptions (a-e) hold in a neighborhood of $\hat x$. However, since the following analysis is local, we can assume without loss of generality that they hold on $\reals^n$. Let $\phi:\reals^n\times\reals \to \reals$ be given by
\[
\phi(x,\gamma) = \Ex\big[\Phi(\bfxi,x,\gamma)\big], ~\mbox{ where } \Phi(\xi,x,\gamma) = \min\big\{0, -\gamma g(\xi,x) - 1 \big\},
\]
which indeed is real-valued because $g(\bfxi,x)$ is integrable by assumption (c).

First, we show that $\phi(\cdot\,,\gamma)$ is differentiable regardless of $\gamma>0$ with the corresponding gradients being continuous jointly in $(x,\gamma)$. Fix $\gamma > 0$ and let $\hat \Xi$ be the set associated with $\hat x$ in assumption (a), i.e., $\prob\{\bfxi\in \hat\Xi\} = 1$ and $g(\xi,\cdot\,)$ is smooth at $\hat x$ for all $\xi \in \hat \Xi$. By \cite[Theorem 7.44]{ShapiroDentchevaRuszczynski.09}, $\phi(\,\cdot\,,\gamma)$ is differentiable at $\hat x$ and
\[
\nabla_x \phi(\hat x, \gamma) = \Ex\big[ \nabla_x \Phi(\bfxi,\hat x,\gamma) \big]
\]
provided that the following two conditions hold:\\ 

{\it Condition 1}: There is a measurable function $\lambda:\reals^m\to [0,\infty)$ and a neighborhood of $\hat x$ such that $\Ex[\lambda(\bfxi)] < \infty$ and
\[
   \big| \Phi(\xi,x,\gamma) - \Phi(\xi,x',\gamma) \big| \leq \lambda(\xi) \|x-x'\|_2
\]
for all $x,x'$ in the neighborhood and all $\xi\in\reals^m$.\\ 

{\it Condition 2}: There is a set $\Xi_1\subset\reals^m$, with $\prob\{\bfxi\in \Xi_1\} = 1$, such that $\Phi(\xi,\cdot,\gamma)$ is differentiable at $\hat x$ for all $\xi\in \Xi_1$.\\

Condition 1 holds with $\lambda(\xi) = \gamma\kappa(\xi)$ by assumption (d). Condition 2 also holds by assumptions (a,e); recall that $g(\xi,\cdot\,)$ is smooth at $\hat x$ for all $\xi \in \hat\Xi$ and $\Phi(\xi,\,\cdot\,,\gamma)$ is smooth at $\hat x$ if
\[
\xi \in \bar\Xi = \hat \Xi \cap \big\{\xi\in\reals^m~\big|~-\gamma g(\xi,\hat x) - 1 \neq 0\big\}.
\]
The event $\bar\Xi$ takes place with probability one and furnishes the required set $\Xi_1$. 

We note that for $\xi\in \bar \Xi$,
\[
\nabla_x \Phi(\xi,\hat x,\gamma) = \begin{cases}
  0 & \mbox{ if } -\gamma g(\xi,\hat x) - 1 > 0\\
  -\gamma \nabla_x g(\xi,\hat x) & \mbox{ if } -\gamma g(\xi,\hat x) - 1 < 0.
\end{cases}
\]

The above arguments hold not only for $\hat x$ but at all other points as well. Thus, $\phi(\,\cdot\,,\gamma)$ is differentiable and the obtained formula for the gradient applies at any point.

By \cite[Theorem 7.43]{ShapiroDentchevaRuszczynski.09}, the mapping $\nabla_x \phi$ is continuous at $(\hat x,\gamma)$ provided that there is a set $\Xi'$, with $\prob\{\bfxi\in \Xi'\} = 1$, such that $\nabla_x \Phi(\xi,\cdot\,,\cdot\,)$ is continuous at $(\hat x, \gamma)$. With $\Xi' = \bar \Xi$, this holds in light of assumption (a). Moreover, the cited theorem requires that there are a measurable function $\eta:\reals^m\to [0,\infty)$, with $\Ex[\eta(\bfxi)] < \infty$, and a set $\Xi''$, with $\prob\{\bfxi\in \Xi''\} = 1$, such that $\|\nabla_x \Phi(\xi,x,\gamma')\|_\infty \leq \eta(\xi)$ for all $\xi\in \Xi''$ and $(x,\gamma')$ in a neighborhood of $(\hat x, \gamma)$. We set $\Xi'' = \bar \Xi$ and utilize $\eta(\xi) = (\gamma + 1)\kappa(\xi)$, which is valid by assumption (d). Thus, we conclude that $\nabla_x \phi$ is continuous at $(\hat x,\gamma)$. Again, we can repeat the argument for other $x$ and $\gamma>0$ and reach the same conclusion.

For the special case of $\gamma = 0$, $\phi(x,0) = -1$ and $\nabla_x \phi(x,0) = 0$. Since $\nabla_x \phi(x^\nu,\gamma^\nu)\to 0$ as $x^\nu\to x$ and $\gamma^\nu\downto 0$, $\nabla_x \phi$ is continuous relative to $\reals^n \times [0,\infty)$.

Following a similar argument, again relying on \cite[Theorem 7.43]{ShapiroDentchevaRuszczynski.09}, we conclude that $\phi$ is continuous.

Second, we consider the function $\psi:\reals^n\to \reals$ given by
\[
\psi(x) = \nmax_{\gamma\geq 0} \phi(x,\gamma).
\]
For each $x$, $\phi(x,\cdot\,)$ is continuous because $\phi$ is continuous. Moreover, $\phi(x,\gamma) \leq (-\gamma \epsilon - 1) \delta$ for all $x\in X$ by assumption (b). Thus, any $\gamma \geq \bar\gamma = (1-\delta)/(\epsilon\delta)$ yields $\phi(x,\gamma) \leq -1$. Since $\phi(x,0) = -1$, this means that the maximization of $\phi(x,\cdot\,)$ over $[0,\infty)$ is attained with a value no greater than $\bar\gamma$ regardless of $x\in X$.  Consequently, $\psi(x)$ is well defined and the maximization can just as well be carried out over the compact set $[0, \bar\gamma]$. We then invoke \cite[Proposition 6.30]{primer} and assumption (f) to conclude that $\psi$ is differentiable at $\hat x$ with $\nabla \psi(\hat x) = \nabla_x \phi(\hat x,\hat\gamma)$. This implies that $-\psi$ is also differentiable at $\hat x$ with gradient
\[
-\nabla_x \phi(\hat x,\hat\gamma) = \Ex\big[ -\hspace{-0.1cm}\nabla_x \Phi(\bfxi,\hat x, \hat\gamma) \big] = \Ex\big[F(\bfxi,\hat x, \hat\gamma) \big].
\]
Since
\[
-\psi(\hat x) = -\nmax_{\gamma\geq 0} \phi(\hat x,\gamma) = \nmin_{\gamma\geq 0} -\phi(\hat x,\gamma) = \bar p(\hat x),
\]
the claim is established.\eop

Interestingly, the proof of the theorem follows a completely different path compared to that of Theorem \ref{thm:discrete}.

Assumption (a) holds if $g(\xi,\cdot\,)$ is smooth for all $\xi\in \reals^m$, but the theorem also allows for certain nonsmooth $g$ of the kind often arising in stochastic programming \cite[Chapter 3]{primer}. For example, if $g(\xi,x) = \max\{0,\beta(x-\xi)\}$ as in the newsvendor problem \cite[Section 1.C]{primer} and $\bfxi$ has a continuous distribution, then assumption (a) still holds. Assumption (b) rules out situations with diminishing probability for positive outcomes of $g(\bfxi,x)$, which is not restrictive in practice. Assumption (e) essentially requires that $g(\bfxi,x)$ has a continuous distribution. As discussed around \eqref{eqn:gamma}, $\hat\gamma$ is ``mostly'' unique so  assumption (f) is rather mild.\\

\state Example 2. Suppose that $\bfxi$ is a uniformly distributed random variable on $[-1,1]$ and $g(\xi,x) = x\xi - 1$. Let $x>1$. Then, the $\alpha$-quantile and $\alpha$-superquantile become $q_\alpha (x) = (2\alpha - 1) x - 1$ and $\bar q_\alpha (x) = \alpha x - 1$. This means that $\bar q_{\bar\alpha}(x) = 0$ implies $\bar \alpha = 1/x$ and then $q_{\bar \alpha}(x) = 1-x$. Moreover, the buffered failure probability $\bar p(x) = 1-1/x$. Certainly, this function is differentiable at any positive $x$ and $\nabla \bar p(x) = 1/x^2$.

Let's now confirm this gradient formula using Theorem \ref{thm:contcase}. By \eqref{eqn:gamma} and the surrounding discussion, $\hat\gamma = 1/(x-1)$ Moreover, one obtains
  \[
F(\xi,x,\hat\gamma) = \begin{cases}
  0 & \mbox{ if } (2-x)/x \geq \xi\\
  \xi/(x-1) & \mbox{ otherwise}
  \end{cases}
  \]
and it only remains to compute the integral
\[
\Ex\big[ F(\bfxi,x,\hat\gamma) \big] = \int_{-1}^1 \tfrac{1}{2}F(\xi,x,\hat\gamma) d\xi = \frac{1}{2}\frac{1}{x-1}\int_{(2-x)/x}^1 \xi\, d\xi = \frac{1}{2}\frac{1}{x-1}   \Big(\frac{1}{2} - \frac{1}{2}\frac{(2-x)^2}{x^2} \Big) = \frac{1}{x^2}.
\]
We have confirmed the formula for $\nabla \bar p(x)$.\eop

\section{Applications}

We consider two applications of Theorem \ref{thm:discrete}. The first one discusses sensitivity analysis of a simple network and illustrates how the formula in \cite{ByunRoyset.21} may not apply; see \eqref{eqn:firstpaper}. The second one presents an optimality condition for constrained minimization of the buffered failure probability.\\

\state Example 3. Consider the network in Figure \ref{fig:fignetwork} consisting of four components. The state of a component $j$ is represented by a Bernoulli random variable $\bfxi_j$ having outcomes 1 and 0, respectively corresponding to its survival and failure, with $\prob\{\bfxi_j = 0\} = \rho\in (0,1/2)$. If a component survives, then it carries a flow equal to its capacity, while it carries no flows if it fails. The capacities of components $1$ and $2$ are $2x$, and those for $3$ and $4$ are $x$, with $x>0$. The hope is that the network delivers at least one unit of flow between the two end nodes. Thus, the quantity of interest becomes the {\em flow-shortfall}: $1$ minus the deliverable flow. We are concerned about the buffered failure probability of the flow-shortfall. With $\xi = (\xi_1, \dots, \xi_4)$, this leads to
\[
g(\xi, x) = \begin{cases}
  1 - 2x & \mbox{ if } \xi_1 = \xi_2 =\xi_3 =\xi_4 =1\\
  1 - x & \mbox{ if } \xi_1 = \xi_2 = \xi_3 = 1, \xi_4 = 0; \mbox{ or } \xi_1 = \xi_2 = \xi_4 = 1, \xi_3 = 0\\
  1 & \mbox{ otherwise.}
\end{cases}
\]
Suppose that $\bfxi = (\bfxi_1, \dots, \bfxi_4)$ and the four random variables are statistically independent. Thus, $\bfxi$ has 16 possible outcomes, denoted by $\xi^1, \dots, \xi^{16}$, but $g(\bfxi,x)$ has only three outcomes $1-2x$, $1-x$ and $1$. This means that the formula \eqref{eqn:firstpaper} from \cite{ByunRoyset.21} does not apply; it remains inapplicable after a consolidation of the probability space into one with only three outcomes.

We consider a particular value for $x$:
\[
\hat x = \frac{1-(1-\rho)^4}{2(1-\rho)^3 \rho},
\]
which is greater than 2 for $\rho \in (0,1/2)$.
For notational convenience, let
\[
\sigma = 1 - (1-\rho)^4 - (1-\rho)^3\rho ~~~\mbox{ and }~~~ \tau = 1 - (1-\rho)^4 - 2(1-\rho)^3\rho.
\]
Via the $\alpha$-superquantile of $g(\bfxi,x)$, we obtain
\[
\bar p(x) = \begin{cases}
2 \sigma x /(2 x - 1) & \mbox{ if } 1/(2(1-\rho)^4) <  x \leq \hat x\\
\tau x /(x - 1) & \mbox{ if } \hat x <  x,
\end{cases}
\]
which is visualized in Figure \ref{fig:barp} for $\rho = 0.1$. Thus, $\bar p$ is neither smooth nor convex nor epi-regular \cite[Defintion 4.57]{primer}. For comparison, Figure \ref{fig:barp} also plots $p(x)$. While not differentiable at $\hat x$, the derivatives of $\bar p$ to the left and right of that point are given by
\[
\nabla \bar p(x) = \begin{cases}
-2\sigma/(2 x - 1)^2  & \mbox{ if } 1/(2(1-\rho)^4) <  x < \hat x\\
-\tau/(x - 1)^2 & \mbox{ if } \hat x <  x.
\end{cases}
\]
The definition of subgradients yields that
\[
\partial \bar p(\hat x) = \bigg\{-\frac{2\sigma}{(2 \hat x - 1)^2} ,~~ -\frac{\tau}{(\hat x - 1)^2}\bigg\}
\]
For $\rho = 0.1$, this produces $\{-0.1073, -0.0392\}$.

Next, we turn to Theorem \ref{thm:discrete} and its (partial) characterization of $\partial \bar p(\hat x)$ in terms of the set
\[
Y(\hat x) = \bigcup_{\hat \gamma\in \Gamma(\hat x)} \Big\{\hat \gamma \sum_{i=1}^{16} p_i \mu_i \nabla_x g(\xi^i, \hat x)~\Big|~ \mu_1, \dots, \mu_{16} \mbox{ satisfy \eqref{eqn:multiplercond} }  \Big\}.
\]
We show that $\partial \bar p(\hat x)$ is a strict subset of $Y(\hat x)$, with the latter being the convex hull of the former.

The distribution function of $g(\bfxi,\hat x)$ is a staircase with a first step of $(1-\rho)^4$ at $1-2\hat x$, a second step of $2(1-\rho)^3 \rho$ at $1-\hat x$ and the remaining step at $1$; see Figure \ref{fig:figcdf}. Since the $(1-\bar p(\hat x))$-quantile of $g(\bfxi,\hat x)$ is $1-2\hat x$ and the value of the distribution function at that point is $1-\bar p(\hat x)$, we conclude that $\Gamma(\hat x)$ is the interval $[1/(2\hat x-1), 1/(\hat x-1)]$; see the discussion around \eqref{eqn:gamma}.

In the following, we order the outcomes of $\bfxi$ as follows: $\xi^1 = (1,1,1,1)$, $\xi^2 = (1,1,0,1)$, $\xi^3 = (1, 1, 1, 0)$ and the remaining outcomes are labeled $\xi^4, \dots, \xi^{16}$ arbitrarily. Then, $g(\xi^1, x) = 1-2 x$, $g(\xi^2, x) = g(\xi^3, x) =  1- x$ and $g(\xi^i, x) = 1$ for $i=4, \dots, 16$.

For any $\hat \gamma \in [1/(2\hat x-1), 1/(\hat x-1)]$, which is a positive number, $\hat \gamma g(\xi^i,\hat x) + 1 > 0$ when $i=4, \dots, 16$. Thus, $\mu_i = 1$ and $\mu_i\nabla_x g(\xi^i,\hat x) = 1 \cdot 0 = 0$ for such $i$.

For $i = 2$ and $3$, $\hat \gamma = 1/(\hat x-1)$ implies that $\hat \gamma g(\xi^i,\hat x) + 1 = 0$. Thus, $\mu_i \in [0,1]$ and $\mu_i\nabla_x g(\xi^i,\hat x) = -\mu_i$ for such $i$ and this value of $\hat \gamma$. For  $\hat \gamma \in [1/(2\hat x-1), 1/(\hat x-1))$, $\hat \gamma g(\xi^i,\hat x) + 1 > 0$, $\mu_i = 1$ and $\mu_i\nabla_x g(\xi^i,\hat x) = 1 \cdot (-1) = -1$.

For $i = 1$, $\hat \gamma = 1/(2\hat x-1)$ implies that $\hat \gamma g(\xi^i,\hat x) + 1 = 0$. Thus, $\mu_i \in [0,1]$ and $\mu_i\nabla_x g(\xi^i,\hat x) = -2 \mu_i$ for such $i$ and this value of $\hat \gamma$. For  $\hat \gamma \in (1/(2\hat x-1), 1/(\hat x-1)]$, $\hat \gamma g(\xi^i,\hat x) + 1 < 0$, $\mu_i = 0$ and $\mu_i\nabla_x g(\xi^i,\hat x) = 0 \cdot (-2) = 0$. Based on the choice of $\hat \gamma$, this leads to three cases.

Case A: $\hat \gamma \in (1/(2\hat x-1), 1/(\hat x-1))$. Then, one has
\[
\sum_{i=1}^{16} p_i \mu_i g(\xi^i, \hat x) = 2(1-\rho)^3 \rho (1-\hat x) + \tau = 0,
\]
which is verified by plugging in the expression for $\hat x$. Hence, the multipliers $\mu_1 = 0$ and $\mu_{i}=1$ for $i>1$ are valid. This leads to the expression
\[
\hat \gamma \sum_{i=1}^{16} p_i \mu_i \nabla_x g(\xi^i, \hat x) = -\hat\gamma 2(1-\rho)^3 \rho.
\]

Case B: $\hat \gamma = 1/(\hat x-1)$. Then, $\mu_1 = 0$ and $\mu_i = 1$ for $i=4, \dots, 16$. The multipliers $\mu_2, \mu_3 \in [0,1]$ and are specified by
\[
\sum_{i=1}^{16} p_i \mu_i g(\xi^i, \hat x) = (\mu_2 + \mu_3)(1-\rho)^3\rho(1-\hat x) + \tau = 0.
\]
This simplifies to $\mu_2 + \mu_3 = 2$ and then $\mu_2 = \mu_3 = 1$. Consequently,
\[
\hat \gamma \sum_{i=1}^{16} p_i \mu_i \nabla_x g(\xi^i, \hat x) = -\frac{2(1-\rho)^3\rho}{\hat x - 1}.
\]

Case C: $\hat \gamma = 1/(2\hat x-1)$. Then, $\mu_i = 1$ for $i>1$ while $\mu_1 \in [0,1]$. The latter is narrowed down by
\[
\sum_{i=1}^{16} p_i \mu_i g(\xi^i, \hat x) = (1-\rho)^4 \mu_1(1-2\hat x) + 2(1-\rho)^3 \rho (1-\hat x) + \tau = 0,
\]
which results in $\mu_1 = 0$. Thus,
\[
\hat \gamma \sum_{i=1}^{16} p_i \mu_i \nabla_x g(\xi^i, \hat x) = - \frac{ 2(1-\rho)^3 \rho}{2\hat x-1}.
\]

Summarizing all the three cases, we see that Theorem \ref{thm:discrete} specifies
\[
Y(\hat x) = \bigg[-\frac{2\sigma}{(2 \hat x - 1)^2} ,~~ -\frac{\tau}{(\hat x - 1)^2}\bigg],
\]
which is the convex hull of $\partial \bar p(\hat x)$.\eop\\

\drawing{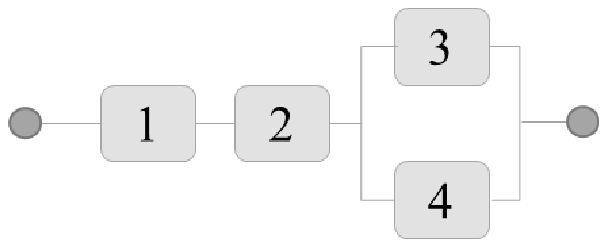}{3.0in} {Network subject to component failure.}{fig:fignetwork}

\drawing{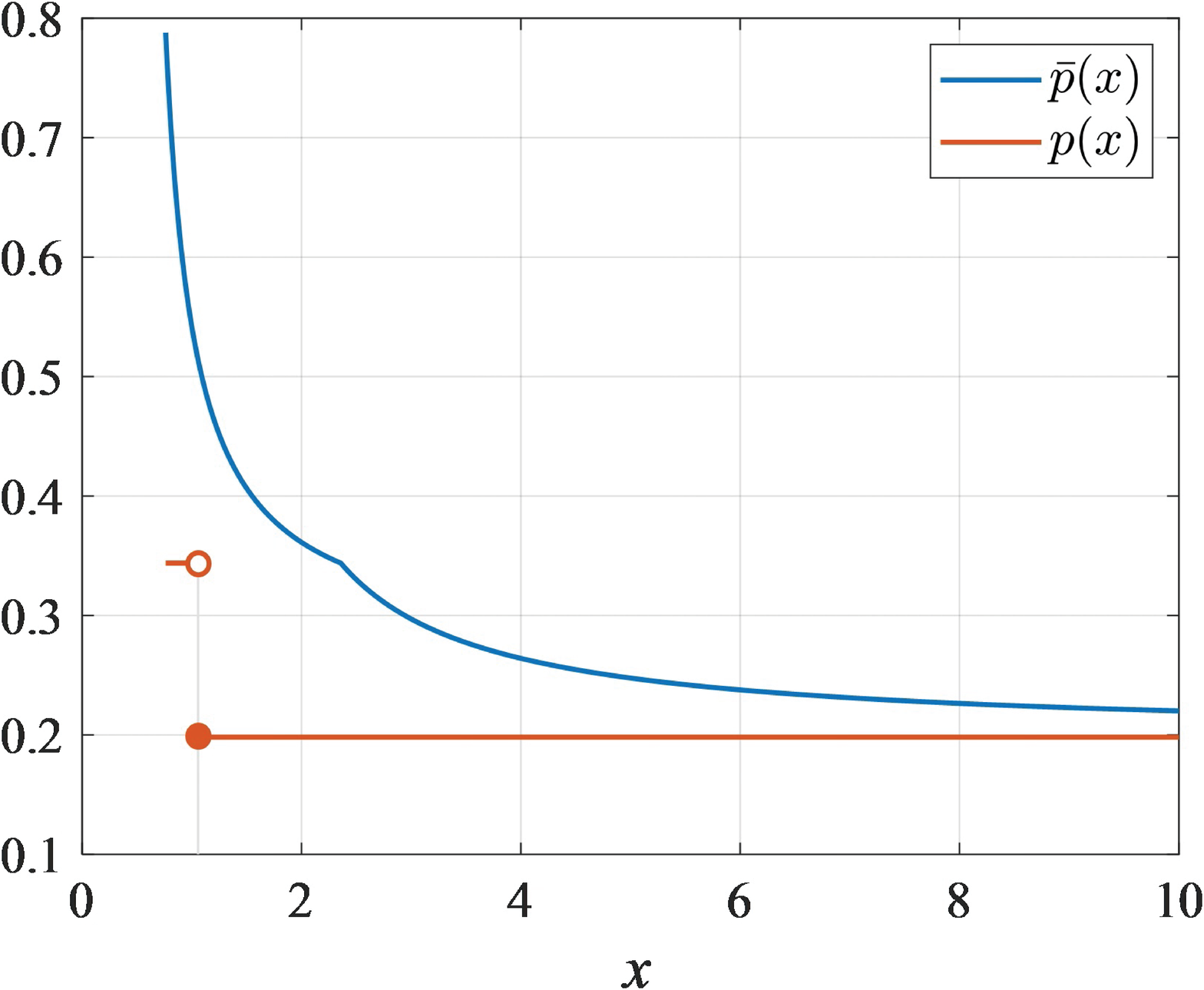}{3.0in} {Buffered failure probability $\bar p(x)$ and failure probability $p(x)$ of flow-shortfall as function of $x$ when $\rho = 0.1$. The kink in the graph occurs at $\hat x = 2.3587$.}{fig:barp}

\drawing{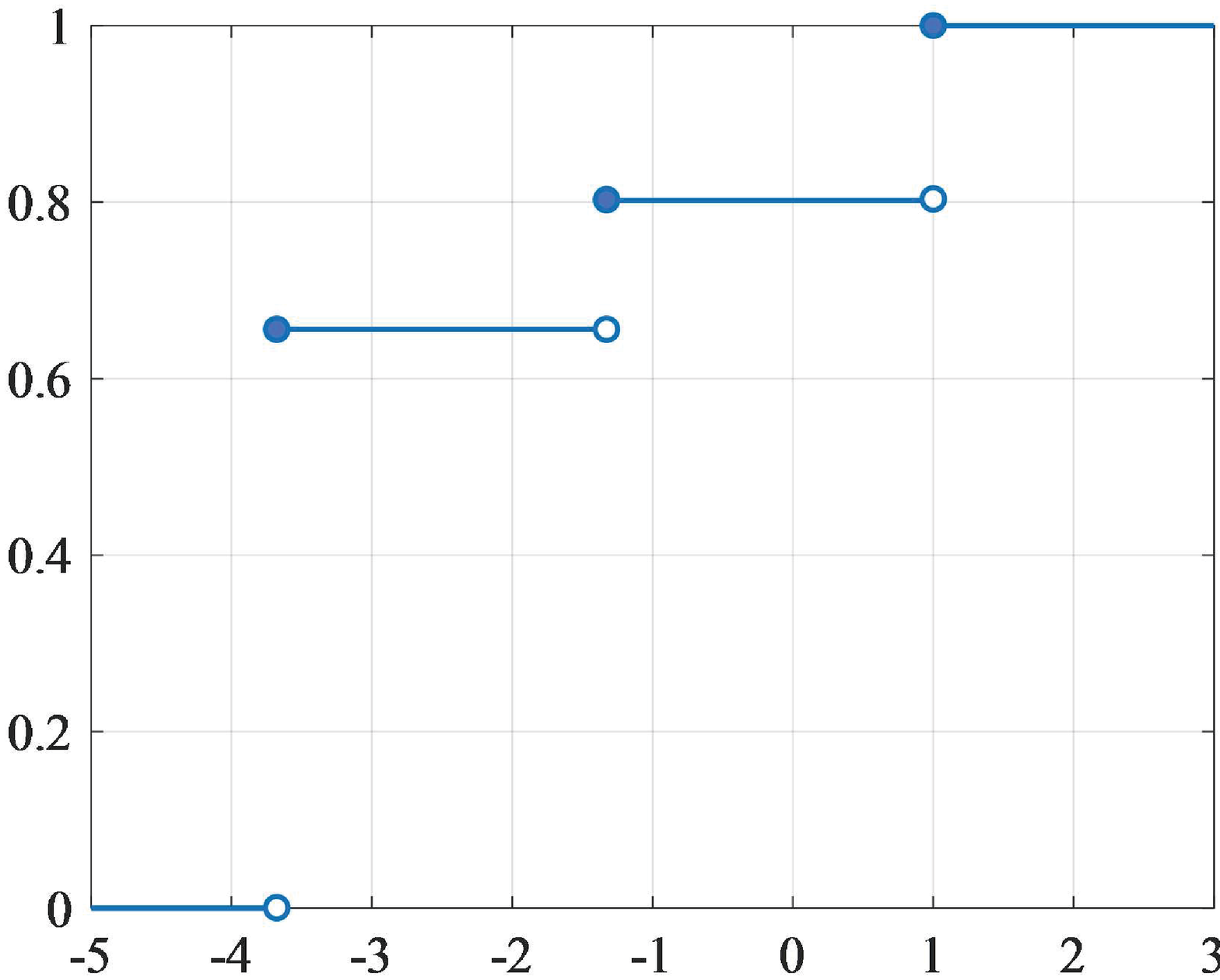}{2.7in} {Distribution function of flow-shortfall $g(\bfxi,\hat x)$ when $\rho = 0.1$.}{fig:figcdf}

\state Example 4. For a nonempty closed set $C \subset\reals^n$, Theorem \ref{thm:discrete} produces an optimality condition for the problem
\[
\nnmin_{x\in C} \bar p(x)
\]
in terms of the normal cone $N_C(x)$ of $C$ at a point $x$ \cite[Definition 4.35]{primer}. (Normal cones are understood in the general or limiting sense; see \cite[Chapter 6]{VaAn}.) 
Specifically, under the assumptions of Theorem \ref{thm:discrete}, if $\hat x$ is a local minimizer of the problem, then there are $\hat\gamma\in \Gamma(\hat x)$ and multipliers $\mu_i \in [0,1], i=1, \dots, \nu$ satisfying \eqref{eqn:multiplercond} such that
\[
-\hat\gamma \sum_{i=1}^\nu p_i \mu_i \nabla_x g(\xi^i,\hat x) \in N_C(\hat x).
\]

This optimality condition holds by the following argument: Let $\hat x$ be a local minimizer of the problem. The Fermat rule \cite[Theorem 4.37]{primer} then asserts that
\[
0 \in \partial (\bar p + \iota_C)(\hat x),
\]
where $\iota_C(x) = 0$ if $x\in C$ and $\iota_C(x) = \infty$ otherwise. Since $\partial \iota_C(\hat x) = N_C(\hat x)$ \cite[Example 4.56]{primer}, the conclusion follows as long as
\begin{equation}\label{eqn:sumrule}
\partial (\bar p + \iota_C)(\hat x) \subset \partial \bar p(\hat x) + N_C(\hat x)
\end{equation}
because then we have the inclusion $0 \in \partial \bar p(\hat x)  + N_C(\hat x)$ and we can invoke Theorem \ref{thm:discrete}.
The inclusion \eqref{eqn:sumrule} holds by a subdifferential sum rule \cite[Proposition 4.67]{primer} provided that the horizon subgradients $\partial^\infty \bar p(\hat x) = \{0\}$ \cite[Definition 4.60]{primer}. This follows by arguing as in the proof of Theorem \ref{thm:discrete}, but now invoking formulas for horizon subgradients instead of for subgradients and employing Theorem 10.6, Corollary 10.9 and Theorem 10.13 in \cite{VaAn}. 

These arguments also confirm that $\bar p$ is locally Lipschitz continuous at $\hat x$; see \cite[Theorem 9.13]{VaAn}.\eop\\

\state Acknowledgement. The research of the first author is supported in part by the Office of Naval Research under MIPR N0001421WX01496 and the Air Force Office of Scientific Research under 18RT0599.

\bibliographystyle{plain}
\bibliography{refs}

\end{document}

%% file: mac.tex
\usepackage{theorem}
\usepackage{enumerate}
\usepackage{array}
\usepackage[reqno]{amsmath}
\usepackage{amssymb}
\usepackage{mathtools}
\usepackage{latexsym}
\usepackage{makeidx}
\usepackage{fancybox}
\input epsf.sty
\usepackage[dvips]{graphicx,color}
\usepackage{showlabels}
\usepackage{subcaption}

\newcommand{\drawing}[4]{
	\begin{figure}[!hbt]
	\begin{center}
	\leavevmode
	\epsfxsize=#2
	\epsfbox{#1}
	\caption{\small #3}
	\label{#4}
	\end{center}
	\end{figure}}

\theoremstyle{change}
\newtheorem{proclaim}{PROCLAIM}[section]
\newtheorem{theorem}[proclaim]{Theorem}

\numberwithin{equation}{section}

\outer\def\proclaim #1. #2\par{\medbreak \noindent{\bf#1.\enspace}{\sl#2}\par
  \ifdim\lastskip<\medskipamount
  \removelastskip\penalty55\medskip\fi}
\def\state #1. { \noindent{\bf#1.\enspace}}
\def\algo #1. { \noindent{\bf#1.\enspace}}

\DeclareMathOperator{\bprob}{b-prob}

\newcommand{\comp}{\,{\raise 1pt \hbox{$\scriptstyle\circ$}}\,}

\newcommand{\reals}{\mathbb{R}}

\newcommand{\natnums}{{{\rm l} \kern -.13em {\rm N} }}

\newcommand{\snats}{{I\kern -.29em N}}
\newcommand{\rats}{{Q\kern -.64em \raise 1pt \hbox{$\scriptstyle |$}\;\,}}
\newcommand{\srats}
	{{Q\kern -.56em \raise 1.2pt \hbox{$\scriptscriptstyle /$}\,}}
\newcommand{\ints}{Z\kern -.46em Z}

\newcommand{\pluss}{\hskip1pt \raise1pt\vbox{\hrule width6pt \vskip1pt \hrule
                    width6pt} \kern-4pt{\lower1pt\hbox{\vrule height6pt
		    \kern1pt\vrule height6pt}}\hskip5pt}
\newcommand{\eop}
	{\hfill{$\vcenter{\hrule height1pt \hbox{\vrule width1pt height5pt
   	 \kern5pt \vrule width1pt} \hrule height1pt}$} \medskip}
\newcommand{\half}
	{{\raisebox{1pt}{$\frac{1}{2}$}}}

\newcommand{\setd}{{ d \kern -.15em l}}
\newcommand{\hatsetd}{ d \hat{\kern -.15em l }}

\renewcommand{\epsilon}{\varepsilon}
\renewcommand{\phi}{\varphi}

\newcommand{\tto}{\;{\lower 1pt \hbox{$\rightarrow$}}\kern -12pt
           \hbox{\raise 2.5pt \hbox{$\rightarrow$}}\;}
\newcommand{\overto}[1]{\,{\raise 0pt\hbox{$\rightarrow$}}\kern -9pt
     \hbox{\lower 3pt \hbox{$\scriptscriptstyle#1$}}\hskip6pt}
\newcommand{\underto}[1]{\,{\lower 1pt\hbox{$\rightarrow$}}\kern -9pt
     \hbox{\raise 4pt \hbox{$\,\scriptscriptstyle#1$}}\hskip7pt}
\newcommand{\bigoverto}[1]{{\raise 0pt\hbox{$\,\longrightarrow$}}\kern -16pt
     \hbox{\lower 3pt \hbox{$\scriptscriptstyle#1$}}\hskip4pt}
\newcommand{\bigunderto}[1]{\,{\lower 1pt\hbox{$\longrightarrow$}}\kern -16pt
     \hbox{\raise 4pt \hbox{$\,\scriptscriptstyle#1$}}\hskip6pt}
\newcommand{\bigbigto}[2]{\,{\raise 0pt\hbox{$\,\longrightarrow$}}\kern -16pt
     \hbox{\lower 3pt \hbox{$\scriptscriptstyle#2$}}\kern -10pt
     \hbox{\raise 4pt \hbox{$\,\scriptscriptstyle#1$}}\hskip7pt}
\newcommand{\downto}{{\raise 1pt \hbox{$\scriptscriptstyle \,\searrow\,$}}}
\newcommand{\upto}{{\raise 1pt \hbox{$\scriptscriptstyle \,\nearrow\,$}}}

\newcommand{\notimply}
	{\quad\hbox{$\Longrightarrow \kern -14pt {/}$}\hskip6pt\quad}

\newcommand{\lto}{\,{\lower 1pt\hbox{$\rightarrow$}}\kern -10pt
     \hbox{\raise 4pt \hbox{$\, \scriptstyle l$}}\hskip7pt}
\newcommand{\eto}{\,{\lower 1pt\hbox{$\rightarrow$}}\kern -11pt
     \hbox{\raise 4pt \hbox{$\, \scriptstyle e$}}\hskip7pt}
\newcommand{\hto}{\,{\lower 1pt\hbox{$\rightarrow$}}\kern -11pt
     \hbox{\raise 4pt \hbox{$\, \scriptstyle h$}}\hskip7pt}
\newcommand{\pto}{\,{\lower 1pt\hbox{$\rightarrow$}}\kern -11pt
     \hbox{\raise 4.5pt \hbox{$\, \scriptstyle p$}}\hskip7pt}
\newcommand{\cto}{\,{\lower 1pt\hbox{$\rightarrow$}}\kern -11pt
     \hbox{\raise 4pt \hbox{$\, \scriptstyle c$}}\hskip7pt}
\newcommand{\gto}{\,{\lower 1pt\hbox{$\rightarrow$}}\kern -11pt
     \hbox{\raise 4.5pt \hbox{$\, \scriptstyle g$}}\hskip7pt}
\newcommand{\sto}{\,{\lower 1pt\hbox{$\rightarrow$}}\kern -11pt
     \hbox{\raise 4pt \hbox{$\, \scriptstyle s$}}\hskip7pt}
\newcommand{\awto}{\,{\lower 1pt\hbox{$\rightarrow$}}\kern -15pt
     \hbox{\raise 4pt \hbox{$\, \scriptstyle aw$}}\hskip7pt}
\def\Nto{\,{\raise 1pt\hbox{$\rightarrow$}}\kern -13pt
     \hbox{\lower 3pt \hbox{$\, \scriptstyle N$}}\hskip7pt}
\def\Cto{\,{\raise 1pt\hbox{$\rightarrow$}}\kern -14pt
     \hbox{\lower 3pt \hbox{$\, \scriptstyle C$}}\hskip7pt}
\def\fto{\,{\raise 1pt\hbox{$\rightarrow$}}\kern -14pt
     \hbox{\lower 3pt \hbox{$\, \scriptstyle f$}}\hskip7pt}

\newcommand{\low}[1]{{\lower1pt \hbox{$\scriptstyle #1$}}}
\newcommand{\loww}[1]{{\lower2pt \hbox{$\scriptstyle #1$}}}
\newcommand{\high}[1]{{\raise1pt \hbox{$\scriptstyle #1$}}}

\newcommand{\ninf}{\mathop{\rm inf}\nolimits}

\newcommand{\nmax}{\mathop{\rm max}\nolimits}
\newcommand{\nmin}{\mathop{\rm min}\nolimits}
\newcommand{\nnmin}{\mathop{\rm minimize}}

\newcommand{\nargmin}{\mathop{\rm argmin}\nolimits}

\newcommand{\prob}{{\mathop{\rm prob}\nolimits}}

\newcommand{\bfxi}{\mbox{\boldmath $\xi$}}

\newcommand{\lwdy}[2]{\mathrel{\mathop
        {\raisebox{0.1ex}{\null$#1$}}{\hbox{\kern -1.0em
	{\raisebox{-0.8ex}{$\scriptstyle{\;\to #2}$}}}}}}
\newcommand{\lwwdy}[2]{\mathrel{\mathop
        {\raisebox{0.2ex}{\null$#1$}}{\hbox{\kern -1.0em
	{\raisebox{-1.1ex}{$\scriptstyle{\;\to #2}$}}}}}}
\newcommand{\slwdy}[2]{\scriptsize{{\mathrel{\mathop
        {\raisebox{0.1ex}{\null$#1$}}{\hbox{\kern -1.0em
	{\raisebox{-0.8ex}{$\scriptstyle{\;\to #2}$}}}}}}}}
\newcommand{\slwwdy}[2]{\scriptsize{{\mathrel{\mathop
        {\raisebox{0.2ex}{\null$#1$}}{\hbox{\kern -1.0em
	{\raisebox{-1.1ex}{$\scriptstyle{\;\to #2}$}}}}}}}}

\definecolor{lightgray}{gray}{0.75}
\definecolor{myred}{rgb}{0.55,0,0}
\definecolor{myblue}{rgb}{0,0,0.5} 
\definecolor{mygreen}{rgb}{0,0.5,0} 
\definecolor{purple}{rgb}{0.5,0,0.5} 
\definecolor{turq}{rgb}{0,0.805,0.816} 
\definecolor{maroon}{rgb}{0.51,0,0}
\definecolor{MAROON}{rgb}{0.51,0,0}
\definecolor{redor}{rgb}{0.78,0.078,0.078}
\definecolor{dgreen}{rgb}{0,0.3,0}

\newcommand{\Ex}{\mathbb{E}}

\newcommand{\bcdot}{\,{\raise .2ex \hbox{$\centerdot$}}\,}

